\documentclass[a4paper,11pt]{article}
\usepackage[active]{srcltx}
\usepackage{amsmath}\usepackage{amssymb}\usepackage{amscd}

\newcommand{\tj}{\tilde{j}}
\newcounter{pic}

\newtheorem{defi}{Definition}\newtheorem{prop}[defi]{Proposition}\newtheorem{lemm}[defi]{Lemma}      
\newtheorem{coro}[defi]{Corollary}   

\newcommand{\ml}[2]{{\mathcal{X}}_{#1}^{#2}}

\newcommand{\xl}[2]{X_{#1}^{#2}}

\newcommand{\lam}{\lambda^{(m)}}
\newcommand{\cc}{c}

\usepackage{empheq}

\begin{document}

$\ $

\begin{center}
\bigskip\bigskip

{\Large\bf On representations of complex reflection groups G(m,1,n)}

\vspace{1.5cm}
{\large {\bf O. V. Ogievetsky$^{\circ\diamond}$\footnote{On leave of absence from P. N. Lebedev Physical Institute, Leninsky Pr. 53,
117924 Moscow, Russia} and L. Poulain d'Andecy$^{\circ}$}}

\vskip 1cm

$\circ\ ${Center of Theoretical Physics, Luminy \\
13288 Marseille, France}

\vspace{.5cm}
$\diamond\ ${J.-V. Poncelet French-Russian Laboratory, UMI 2615 du CNRS, Independent University of Moscow, 11 B. Vlasievski per., 
119002 Moscow, Russia}
\end{center}

\vskip 1cm
\begin{abstract}
An inductive approach to the representation theory of the chain of the complex reflection groups $G(m,1,n)$
is presented. We obtain the Jucys-Murphy elements of $G(m,1,n)$ from the Jucys--Murphy elements of the cyclotomic Hecke algebra, and study their common spectrum using representations of a degenerate cyclotomic affine Hecke algebra. Representations of $G(m,1,n)$ are constructed with the help of a new associative algebra whose underlying vector space is the tensor product of 
the group ring $\mathbb{C}G(m,1,n)$ with a free associative algebra generated by the standard $m$-tableaux.
\end{abstract}

\section{{\hspace{-0.55cm}.\hspace{0.55cm}}Introduction\vspace{.25cm}}

The complex reflection groups generalize the 
Coxeter groups and the complete list of irreducible finite complex reflection groups consists of the series of groups denoted $G(m,p,n)$, where 
$m,p,n$ are positive integers such that $p$ divides $m$, and $34$ exceptional groups \cite{ST}. Analogues of the Hecke algebra and the braid group exist for all complex reflection groups.

\vskip .2cm
The Hecke algebra of $G(m,1,n)$, which we denote by $H(m,1,n)$, has been introduced 
in \cite{AK,Bro-M,C2} and is called the cyclotomic Hecke algebra. For $m=1$ this is the Hecke algebra of type $A$ and for $m=2$ this is the Hecke 
algebra of type $B$. The representation theory of the algebras $H(m,1,n)$ was developed in \cite{AK} 
(and in \cite{H} for the Hecke algebra of type $B$); in \cite{OPdA}, an inductive approach, \`a la Okounkov--Vershik \cite{OV}, to the representation theory of the chain, with respect to $n$, of the algebras $H(m,1,n)$ was suggested.

\vskip .2cm
In this paper, we present an inductive approach, in the spirit of \cite{OV}, to the representation theory of the chain of the complex reflection groups $G(m,1,n)$. We construct representations of $G(m,1,n)$ with the help of a certain associative algebra, denoted further by $\mathfrak{T}$; the underlying vector space of $\mathfrak{T}$ is
the tensor product of the group ring $\mathbb{C}G(m,1,n)$ with the free algebra generated by the standard $m$-tableaux of shape $\lambda^{(m)}$ for all $m$-partitions $\lambda^{(m)}$ of length $n$. 

\vskip .2cm
The representation theory of the groups $G(m,1,n)$ is well known, see, \cite{Sp} or, \emph{e.g.}, \cite{Mac} (for the Coxeter groups of type $B$, 
that is, for the groups $G(2,1,n)$, the representation theory was developed in \cite{Y}). Also, the representation theory of $G(m,1,n)$ can be directly deduced from the representation theory of $H(m,1,n)$ by taking the following limit
\begin{equation}\label{clali}
\left\{\begin{array}{cl}
v_i\to\xi_i & \textrm{for $i=1,\dots,m$, where the $\xi_i$ are distinct $m^{\textrm{th}}$ roots of unity\ ,}\\[.5em]
q\to\pm1\ ,&
 \end{array}\right.
\end{equation}
 in the formulas \cite{OPdA} for the matrix elements of the generators;
the parameters $q, v_1,\dots,v_m$ enter the definition of the cyclotomic Hecke algebra in the following way: the algebra $H(m,1,n)$ is generated by $\tau,\sigma_1,\dots,\sigma_{n-1}$ with the defining relations  
\begin{equation}\label{pres-hec-cyc}
\left\{\begin{array}{ll}
 \sigma_i\sigma_{i+1}\sigma_i=\sigma_{i+1}\sigma_i\sigma_{i+1} & \text{for $i=1,\dots,n-2$}\ ,\\[0.3em]
 \sigma_i\sigma_j=\sigma_j\sigma_i & \text{for $i,j=1,\dots,n-1$ such that $|i-j|>1$}\ ,\\[0.3em]
 \sigma_i^2=(q-q^{-1})\sigma_i+1 & \text{for $i=1,\dots,n-1$}\ ,\\[0.3em]
 \tau\sigma_1\tau\sigma_1=\sigma_1\tau\sigma_1\tau\ ,\\[0.3em]
 \tau\sigma_i=\sigma_i\tau & \text{for $i>1$}\ ,\\[0.3em]
 (\tau-v_1)\dots(\tau-v_m)=0\ .
\end{array}\right.
\end{equation}

However it is interesting to take the ``classical limit" of the whole approach, developed in \cite{OPdA}, establishing thereby the inductive approach to the representation theory of the groups $G(m,1,n)$.
We insist here on the connection of the treatment for the groups $G(m,1,n)$ with the treatment for the algebras $H(m,1,n)$. Also the existence of the algebra $\mathfrak{T}$ is of independent interest and the construction of the representations presented here appears to be new. 

\vskip .2cm
The representation theory of a more general class of groups, namely,  
the wreath products of finite groups by the symmetric groups, was built, in the spirit of \cite{OV}, in  \cite{Pus}. The construction in \cite{Pus} is worked out 
within the group theory. In this paper, we shall see how this approach
is restored -- on the example of the groups $G(m,1,n)$, the wreath products of the cyclic groups by the symmetric groups -- 
in the classical limit of the construction developed in \cite{OPdA} for $H(m,1,n)$.
We shall see that there are certain subtleties in passing to the classical ``group" situation (one should be careful about the order of taking limits {\it etc}). As it often happens the classical situation is more complicated than the quantum one. 

\vskip .2cm
The branching rules for the chain of wreath products $A\wr S_n$, where $A$ is a finite group, are multiplicity free if and only if the group $A$ is abelian. The chain of groups $G(m,1,n)$ provides a simplest example of a multiplicity free chain of wreath products $A\wr S_n$. 

\vskip .2cm
First we present the classical Jucys--Murphy elements of the 
group ring of $G(m,1,n)$, which we obtain as classical limits of certain expressions involving the Jucys--Murphy elements of $H(m,1,n)$ (a similar process was used in 
\cite{Ram-JM} for the Weyl groups). The Jucys--Murphy elements were defined in \cite{Pus,Wang} 
for the wreath product of any finite group $A$ by the symmetric group. 
The Jucys--Murphy elements obtained by a limiting procedure coincide with those from \cite{Pus,Wang} if we choose $A$ to be the cyclic group. As in the non-degenerate situation, the usage of the Jucys--Murphy elements is the main tool in our construction of the representation theory.

\vskip .2cm
The representation 
theory of the Hecke algebras $H(m,1,n)$ in the inductive approach requires, for all $m$, the study of representations of the same affine Hecke algebra of type A. However, in 
the classical limit, a certain version of degenerate affine cyclotomic Hecke algebras, which we denote ${\mathfrak{A}}_{m,n}$ in the text, appears as the classical counterpart of the affine Hecke algebra; the 
representations of the simplest degenerate affine cyclotomic Hecke algebra ${\mathfrak{A}}_{m,2}$ serve for the study of the representation theory of $G(m,1,n)$ -- the classical limit of $H(m,1,n)$. 

\vskip .2cm
The algebras ${\mathfrak{A}}_{m,n}$ for all $m=1,2,\dots$ can be obtained by a certain limiting procedure (the details are below)
from one and the same affine Hecke algebra $\hat{H}_n$. 

\vskip .2cm
We show that the Jucys--Murphy elements of the group ring of $G(m,1,n)$ are images of the ``universal" Jucys--Murphy elements living 
in the algebra ${\mathfrak{A}}_{m,n}$. 
We verify, on the classical level, the commutativity of the double set 
$\{ x_1,\tilde{x}_1,\dots, x_n,\tilde{x}_n\}$ of elements in the algebra $\mathfrak{A}_{m,n}$. The algebra ${\mathfrak{A}}_{m,n}$ turns out to coincide with a particular case of the wreath Hecke algebra defined in \cite{WaWa}, see also \cite{Ram-Sh}. 
We do not include the commutativity of the set $\{ x_1,\tilde{x}_1,\dots, x_n,\tilde{x}_n\}$ of elements in the defining relations of $\mathfrak{A}_{m,n}$, contrary to the definition of the wreath Hecke algebra of \cite{WaWa}; as a corollary of the results here, the two algebras are in fact isomorphic.

\vskip .2cm
The representation theory of the 
simplest non-trivial degenerate cyclotomic affine Hecke algebra, the algebra ${\mathfrak{A}}_{m,2}$, carries an important information about the recursive properties of the Jucys--Murphy elements of the group ring of $G(m,1,n)$.
We present the list of irreducible representations with diagonalizable $x_1,\tilde{x}_1,x_2$ and $\tilde{x}_2$ of the algebra $\mathfrak{A}_{m,2}$, and then use it to study the spectrum of the Jucys--Murphy elements. In particular, we characterize the sets of common eigenvalues of the 
Jucys--Murphy elements of the group ring of $G(m,1,n)$ and establish then the relation with the $m$-partitions and the $m$-tableaux. Once the representations of the group $G(m,1,n)$ are 
constructed, we verify by counting dimensions that all irreducible representations are obtained within this approach.

\vskip .2cm
The representations of $G(m,1,n)$ constructed with the help of the algebra $\mathfrak{T}$ (the tensor product of the algebra $\mathbb{C}G(m,1,n)$ with a free associative algebra generated by the standard 
$m$-tableaux corresponding to $m$-partitions of $n$) are analogues, for $G(m,1,n)$, of the semi-normal representations of the symmetric group. 
We determine the $G(m,1,n)$-invariant Hermitian scalar product on the representation spaces
and describe analogues of the orthogonal representations of the symmetric group.

\vskip .2cm
We briefly outline the organization of the paper. In Section \ref{subsec-cla1}, we recall the standard presentation of the complex reflection groups $G(m,1,n)$, and 
define, by a limiting procedure, the Jucys--Murphy elements of $\mathbb{C}G(m,1,n)$. In Section \ref{subsec-cla3}, we present the degenerate cyclotomic affine Hecke algebra ${\mathfrak{A}}_{m,n}$ and show that the Jucys--Murphy elements of the group ring of $G(m,1,n)$ are images of the ``universal" Jucys--Murphy elements living in ${\mathfrak{A}}_{m,n}$. In Section \ref{sec-spec}, we study the spectrum of the Jucys--Murphy elements of the group ring of $G(m,1,n)$ and establish the connection with the standard $m$-tableaux. The representations of the group $G(m,1,n)$ are 
constructed in Section \ref{suserep-c} with the help of the algebra $\mathfrak{T}$. We conclude our study of the representation theory of $G(m,1,n)$ in Section \ref{subsec-comp-c}, in which we give the completeness result and some consequences. In the Appendix
we study the intertwining operators (introduced in \cite{WaWa}) in the degenerate cyclotomic affine Hecke algebra. The intertwining operators provide certain information about the spectrum of the Jucys--Murphy elements. We 
explain how to obtain these intertwining operators by taking the classical limit of appropriate intertwining operators of the non-degenerate affine Hecke algebra.

\vskip .2cm 
We refer to \cite{OPdA2} for the definitions and results concerning the cyclotomic Hecke algebra $H(m,1,n)$ invoked here, as well as
for the proofs, skipped in the present text.

\paragraph{Notation.} In this article, the ground field is the field $\mathbb{C}$ of complex numbers. 
The spectrum of an operator ${\cal{T}}$ is denoted by ${\mathrm{Spec}}({\cal{T}})$.
We denote, for two integers $k,l\in\mathbb{Z}$ with $k<l$, by  $[k,l]$ the set of integers $\{k,k+1,\dots,l-1,l\}$.

\section{{\hspace{-0.55cm}.\hspace{0.55cm}}Complex reflection group $G(m,1,n)$ and Jucys--Murphy elements}\label{subsec-cla1}

The group $G(m,1,n)$ is generated by the elements $t$, $s_1$, $\dots$, $s_{n-1}$ with the defining relations:
\begin{equation}\label{def2a}
\left\{\begin{array}{ll}
s_is_{i+1}s_i=s_{i+1}s_is_{i+1} & \textrm{for $i=1,\dots,n-2$\ ,}\\[.2em]
s_is_j=s_js_i & \textrm{for $i,j=1,\dots,n-1$ such that $|i-j|>1$}\ ,\\[.2em]
s_i^2=1 & \textrm{for $i=1,\dots,n-1$}
 \end{array}\right.
\end{equation}
and
\begin{equation}\label{def2b}
\left\{\begin{array}{ll}
ts_1ts_1=s_1ts_1t\ ,\\[.2em]
ts_i=s_it & \textrm{for $i>1$\ ,}\\[.2em]
t^m=1\ . & 
\end{array}\right.
\end{equation}

The group $G(m,1,n)$ is isomorphic to the group $C_m\wr S_n$, the wreath product of the cyclic group with $m$ elements, $C_m$, by the symmetric group $S_n$.
; its order is equal to $m^nn!$.
The subgroup of the group $G(m,1,n)$ generated by the elements $s_1$, $\dots$, $s_{n-1}$ is isomorphic to the symmetric group $S_n$. 

\paragraph{Jucys--Murphy elements.} In the following we identify the generators $\tau$, $\sigma_1$, $\dots$, $\sigma_{n-1}$ of $H(m,1,n)$
with respectively $t$, $s_1$ ,$\dots$, $s_{n-1}$ as soon as we have taken the classical limit (\ref{clali}).

Recall that the Jucys--Murphy elements $J_i$ of $H(m,1,n)$ are defined by the initial condition $J_1=\tau$ and the recursion $J_{i+1}=\sigma_iJ_i\sigma_i$, $i=1,\dots,n-1$. We define the following classical analogues of the elements $J_i$: 
\begin{equation}\label{JM-Bn1}
j_i:=\lim_{q\to1}\lim_{v_i\to\xi_i}\bigl(J_i\bigr)\ ,
\end{equation}
and
\begin{equation}\label{JM-Bn2}
\tj_i:=\frac{1}{m}\ \lim_{q\to1}\lim_{v_i\to\xi_i}\left(\frac{J_i^m-1}{q-q^{-1}}\right)\ .
\end{equation}
Attention: the order of taking limits here is important, we first take the limit with respect to the variables $v_i$ and then with respect to $q$; it is maybe
more instructive to write (\ref{JM-Bn2}) in the form $\tj_i:=\frac{1}{m}\ \lim\limits_{q\to1}\ \frac{\lim\limits_{v_i\to\xi_i}\left( J_i^m-1\right)}{q-q^{-1}}\ $. 

\section{{\hspace{-0.55cm}.\hspace{0.55cm}}Degenerate cyclotomic affine Hecke algebra}\label{subsec-cla3}

The Jucys--Murphy elements of the cyclotomic Hecke algebra $H(m,1,n)$ are the images of the ``universal" Jucys--Murphy elements of the affine Hecke algebra. Similarly, the elements $j_i$ and $\tj_i$ are the images of certain elements of a ``universal" degenerate cyclotomic affine Hecke
algebra, which we introduce here.

\begin{defi}
{\hspace{-.2cm}.\hspace{.2cm}}
 \label{def-deg-affHec}
Let $\mathfrak{A}_{m,n}$ be the algebra generated by ${\overline s}_1,\dots,{\overline s}_{n-1}$ and two more generators, 
$x_1$ and $\tilde{x}_1$; the defining relations we introduce in three steps. 
First, there are defining relations, involving the generators ${\overline s}_1,\dots,{\overline s}_{n-1}$ only: 
\begin{equation}\label{deg-affHec-na}
\left\{\begin{array}{ll}
{\overline s}_i{\overline s}_{i+1}{\overline s}_i={\overline s}_{i+1}{\overline s}_i{\overline s}_{i+1} & \textrm{for $i=1,\dots,n-2$\ ,}\\[0.4em]
{\overline s}_i{\overline s}_j={\overline s}_j{\overline s}_i & \textrm{for $i,j=1,\dots,n-1$ such that $|i-j|>1$\ ,}\\[0.4em]
{\overline s}_i^2=1 & \textrm{for $i=1,\dots,n-1$\ ;}
\end{array}\right.
\end{equation}
second, there are relations concerning the addition of the generator $x_1$:
\begin{equation}\label{deg-affHec-nb}
\left\{\begin{array}{l}
x_1{\overline s}_1x_1{\overline s}_1={\overline s}_1x_1{\overline s}_1x_1\ ,\\[0.4em]
x_1{\overline s}_i={\overline s}_ix_1\ \quad\textrm{for $i>1$\ ,}\\[0.4em]
x_1^m=1\ ;
\end{array}\right.
\end{equation}
the third group of relations concerns the addition of the last generator $\tilde{x}_1$:
\begin{equation}\label{deg-affHec-nc}
\left\{\begin{array}{l}
 \tilde{x}_1({\overline s}_1\tilde{x}_1{\overline s}_1+\frac{1}{m}\sum\limits_{p=1}^mx_1^p{\overline s}_1x_1^{-p})=
({\overline s}_1\tilde{x}_1{\overline s}_1+\frac{1}{m}\sum\limits_{p=1}^mx_1^p{\overline s}_1x_1^{-p})\tilde{x}_1\ ,\\[0.6em]
\tilde{x}_1{\overline s}_i={\overline s}_i\tilde{x}_1\ \quad\textrm{for $i>1$\ ,}\\[0.4em]
\tilde{x}_1x_1=x_1\tilde{x}_1\ ,\\[0.4em]
\tilde{x}_1{\overline s}_1x_1{\overline s}_1={\overline s}_1x_1{\overline s}_1\tilde{x}_1\ .
\end{array}\right.
\end{equation}

\vskip .1cm\noindent
We call the algebra $\mathfrak{A}_{m,n}$ the degenerate cyclotomic affine Hecke algebra.
\end{defi}

Due to the relations (\ref{deg-affHec-na})--(\ref{deg-affHec-nb}) there is a homomorphism 
\begin{equation}\label{cyctdegaf}\hat{\iota}\, :\ \mathbb{C}G(m,1,n)\to \mathfrak{A}_{m,n}\ ,\quad \hat{\iota}(s_i)={\overline s}_i
\ \ {\text{for}}\ \ i=1,\dots,n-1\ ,\ 
\hat{\iota}(t)=x_1\ .\end{equation}
Let $\pi$ be a map from the set of generators $\{ {\overline s}_1,\dots,{\overline s}_{n-1},x_1,\tilde{x}_1\}$ to $\mathbb{C}G(m,1,n)$
defined by
\begin{equation}\label{degaftcycl} \pi\, :\ {\overline s}_i\mapsto s_i\ \ {\text{for}}\ \ i=1,\dots,n-1\ ,\ x_1\mapsto t\ ,\ \tilde{x}_1\mapsto 0\ . \end{equation}
Clearly, $\pi$ extends to a homomorphism, which we denote by the same symbol $\pi$, from the algebra $\mathfrak{A}_{m,n}$ to $\mathbb{C}G(m,1,n)$  
(the homomorphism $\pi$ is well defined since the relations (\ref{deg-affHec-nc}) are trivially satisfied when one sends $\tilde{x}_1$ to 0).
Moreover, the composition $\pi\circ\hat{\iota}$ leaves the generators of $G(m,1,n)$ invariant and is therefore the identity homomorphism of the algebra 
$\mathbb{C}G(m,1,n)$; in particular, the map $\hat{\iota}$ is injective or, equivalently,
the subalgebra of $\mathfrak{A}_{m,n}$ generated by the elements ${\overline s}_1,\dots,{\overline s}_{n-1}$ and $x_1$ is isomorphic to the algebra $\mathbb{C}G(m,1,n)$. 

\vskip .2cm
Define ``higher" elements $x_i$ and $\tilde{x}_i$ for $i=2,\dots,n$ by
\begin{equation}\label{hixk}x_{i+1}={\overline s}_ix_i{\overline s}_i\, \quad\textrm{and}\quad
\tilde{x}_{i+1}={\overline s}_i\tilde{x}_i{\overline s}_i+\frac{1}{m}\sum\limits_{p=1}^mx_i^p{\overline s}_ix_i^{-p}\, ,\ \ i=1,\dots,n-1\ .\end{equation}
The first relation in (\ref{deg-affHec-nb}) and the first and fourth relations in (\ref{deg-affHec-nc}) can be rewritten, respectively, as
\begin{equation}\label{x1cox2}x_1x_2=x_2x_1\ ,\ \ \tilde{x}_1\tilde{x}_2=\tilde{x}_2\tilde{x}_1\ \ {\text{and}}\ \ \tilde{x}_1x_2=x_2\tilde{x}_1\ .\end{equation}

\begin{prop}
{\hspace{-.2cm}.\hspace{.2cm}}
 \label{lem1}
We have
\begin{equation}\label{jmicyah}\pi (x_i)=j_i\ \ {\text{and}}\ \ \pi (\tilde{x}_i)=\tilde{j}_i\ .\end{equation}
\end{prop}

Since the Jucys--Murphy elements $J_i$ 
commute in the algebra $H(m,1,n)$, it follows from the definitions (\ref{JM-Bn1}) and (\ref{JM-Bn2}) that the elements 
$j_i$, $i=1,\dots,n$, and the elements $\tj_i$, $i=1,\dots,n$, form together a commutative set. We did not include the commutativity of the corresponding 
set, formed by the elements $x_i$, $i=1,\dots,n$, and the elements $\tilde{x}_i$, $i=1,\dots,n$, in the 
defining 
relations for the algebra  $\mathfrak{A}_{m,n}$: 
the commutativity of this set (and therefore, by the Lemma \ref{lem1}, of its image under the morphism $\pi$, that is,
of the set formed by the elements $j_i$, $i=1,\dots,n$, and $\tj_i$, $i=1,\dots,n$) 
follows, as we shall now see, from the relations (\ref{deg-affHec-na})--(\ref{deg-affHec-nc}).

\begin{lemm}
{\hspace{-.2cm}.\hspace{.2cm}} \label{lem2}
The relations (\ref{deg-affHec-na})--(\ref{deg-affHec-nc}) imply that $x_i$ and $\tilde{x}_i$ commute with ${\overline s}_k$ for $k>i$ and $k<i-1$.
\end{lemm}

\begin{prop}
{\hspace{-.2cm}.\hspace{.2cm}} \label{prop5}
The relations (\ref{deg-affHec-na})--(\ref{deg-affHec-nc}) imply that:
\begin{equation}\label{corexyfdr}x_kx_l=x_lx_k\ ,\ \tilde{x}_k\tilde{x}_l=\tilde{x}_l\tilde{x}_k\ \ \textrm{and}\ \ x_k\tilde{x}_l=
\tilde{x}_lx_k\ \ \textrm{for all $\ k,l=1,\dots,n$\ .}\end{equation}
\end{prop}

\section{{\hspace{-0.55cm}.\hspace{0.55cm}}Spectrum of the Jucys--Murphy elements}\label{sec-spec} 

\subsection{{\hspace{-0.50cm}.\hspace{0.50cm}}Representations of algebra $\mathfrak{A}_{m,2}$}\label{subsec-cla4} 

The important step in the understanding of the spectrum of the Jucys--Murphy elements and in the 
construction of representations is the analysis of the representations of the smallest non-trivial degenerate cyclotomic affine Hecke algebra, 
the algebra $\mathfrak{A}_{m,2}$.
Here we present the list of irreducible representations with diagonalizable $x_1,\tilde{x}_1,x_2$ and $\tilde{x}_2$ of the algebra $\mathfrak{A}_{m,2}$.

\bigskip
Consider the algebra $\mathfrak{A}_{m,2}$ generated by $x$, $y$, $\tilde{x}$, $\tilde{y}$ and $s$ with the relations:
\begin{equation}\label{deg-affHec}\left\{\begin{array}{l}xy=yx\ ,\quad\tilde{x}\tilde{y}=\tilde{y}\tilde{x}\ ,\quad x\tilde{x}=\tilde{x}x\ ,
\quad y\tilde{x}=\tilde{x}y\ ,
\\[.5em]
y=sxs\ ,\quad x^m=1\ ,\quad\tilde{y}=s\tilde{x}s+\frac{1}{m}\sum\limits_{p=1}^mx^psx^{m-p}\ ,\quad s^2=1\ .
\end{array}\right.\end{equation} 
For all $i=1,\dots,n-1$, the subalgebra of $\mathbb{C}G(m,1,n)$ generated by $j_i$, $j_{i+1}$, $\tj_i$, $\tj_{i+1}$ and $s_i$ is a quotient of the 
algebra $\mathfrak{A}_{m,2}$. For $m=1$ the algebra $\mathfrak{A}_{m,2}$ reduces to the degenerate affine Hecke algebra studied in \cite{OV} for the 
representation theory of the symmetric groups $S_n$.

\vskip .2cm
The four elements $x$, $y$, $\tilde{x}$ and $\tilde{y}$ pairwise commute, see the Proposition \ref{prop5}. We investigate irreducible 
representations of the algebra $\mathfrak{A}_{m,2}$ with diagonalizable $x$, $y$, $\tilde{x}$ and $\tilde{y}$. We obtain the following classification.
\begin{itemize}
 \item One-dimensional representations. The action of  the generators is given by
\begin{equation}\label{mat-d1'}
x\mapsto a\ ,\quad y\mapsto a\ ,\quad \tilde{x}\mapsto\tilde{a}\ ,\quad \tilde{y}\mapsto\tilde{a}+\epsilon\ ,\quad s\mapsto\epsilon\ ,\end{equation}
where $a^m=1$ and $\epsilon^2=1$.
\item Two kinds of two-dimensional representations. In the first one, the matrices of the generators of the algebra $\mathfrak{A}_{m,2}$ are given by 
\begin{equation}\label{mat-d2'}s\mapsto\left(\begin{array}{cc}0 & 1\\ 1 & 0\end{array}\right), \  
x\mapsto\left(\begin{array}{cc}a & 0\\ 0 & b\end{array}\right), \  y\mapsto\left(\begin{array}{cc}b & 0\\ 0 & a\end{array}\right), \ \tilde{x}\mapsto\left(\begin{array}{cc}\tilde{a} & 0\\ 0 & \tilde{b}\end{array}\right), \  
\tilde{y}\mapsto\left(\begin{array}{cc}\tilde{b} & 0\\ 0 & \tilde{a}\end{array}\right)\,,\end{equation}
where $a^m=b^m=1$ and $a\neq b$.
In the second one, the matrices of the generators of the algebra $\mathfrak{A}_{m,2}$ are given by 
\begin{equation}\label{mat-d2''}s\mapsto\left(\begin{array}{cc}
(\tilde{b}-\tilde{a})^{-1} & 1-(\tilde{b}-\tilde{a})^{-2}\\[.1em] 1 & -(\tilde{b}-\tilde{a})^{-1}\end{array}\right),\  
x,y\mapsto\left(\begin{array}{cc}a & 0\\ 0 & a\end{array}\right),\ 
\tilde{x}\mapsto\left(\begin{array}{cc}\tilde{a} & 0\\ 0 & \tilde{b}\end{array}\right),\ 
\tilde{y}\mapsto\left(\begin{array}{cc}\tilde{b} & 0\\ 0 & \tilde{a}\end{array}\right)\,,\end{equation}
where $a^m=1$ and $\tilde{b}\neq\tilde{a}$. The representation (\ref{mat-d2''}) is irreducible if and only if  $\tilde{b}\neq\tilde{a}\pm1$.
\end{itemize}

\subsection{{\hspace{-0.50cm}.\hspace{0.50cm}}Classical spectrum}\label{subsec-cla5} 

We begin by constructing representations of $\mathbb{C}G(m,1,n)$ verifying two conditions. First, the classical 
Jucys--Murphy elements  $j_1,\dots,j_n,\tj_1,\dots,\tj_n$ are represented by semi-simple (diagonalizable) operators. Second, for every $i=1,\dots,n-1$ 
the action of the subalgebra generated by $j_i$, $j_{i+1}$, $\tj_i$, $\tj_{i+1}$ and $s_i$ is completely reducible. We shall use the name $C$-representations ($C$ is the first letter in ``completely reducible") for these representations.  At the end of the construction we shall see that all irreducible representations of $\mathbb{C}G(m,1,n)$ are $C$-representations.

\vskip .2cm
We denote ${\mathrm{Spec}}
\left(\begin{array}{ccc}j_1&,\, \dots\, ,&j_n\\[.5em]
\tj_1&,\, \dots\, ,&\tj_n\end{array}\right)$ the set of common eigenvalues of the elements 
$j_1,\tj_1, \dots,j_n,\tj_n$ in the $C$-representations:
\begin{equation}\label{2byns}\Lambda=\left(\begin{array}{ccc}a^{(\Lambda)}_1&,\, \dots\, ,&a^{(\Lambda)}_n\\[.5em]
\tilde{a}^{(\Lambda)}_1&,\, \dots\, ,&\tilde{a}^{(\Lambda)}_n\end{array}\right)\end{equation}
belongs to ${\mathrm{Spec}}
\left(\begin{array}{ccc}j_1&,\, \dots\, ,&j_n\\[.5em]
\tj_1&,\, \dots\, ,&\tj_n\end{array}\right)$ if there is a vector $e_{\Lambda}$ in the space of a $C$-representation 
such that $j_i(e_{\Lambda})=a^{(\Lambda)}_ie_{\Lambda}$ and $\tj_i(e_{\Lambda})=\tilde{a}^{(\Lambda)}_ie_{\Lambda}$ for all $i=1,\dots,n$. 

\vskip .2cm
The elements $j_i$ and $\tj_i$ commute with $s_k$ for $k>i$ and $k<i-1$ (see the Lemma \ref{lem2}) 
which implies that the action of $s_k$ on ${\mathrm{Spec}}
\left(\begin{array}{ccc}j_1&,\, \dots\, ,&j_n\\[.5em]
\tj_1&,\, \dots\, ,&\tj_n\end{array}\right)$ is ``local" in the sense that $s_k(e_{\Lambda})$ is a linear combination of $e_{\Lambda'}$ such that 
${{a^{(\Lambda')}_i=a^{(\Lambda)}_i}^{\phantom{A}}}^{\phantom{A}}$ 
\hspace{-0.55cm} and $\tilde{a}^{(\Lambda')}_i=\tilde{a}^{(\Lambda)}_i$ for $i\neq k,k+1$. 

\vskip .2cm
The $2\times n$ arrays
(\ref{2byns}) we shall call strings, keeping the name ``string" used for a set of common eigenvalues of the Jucys--Murphy elements for the algebra $H(m,1,n)$.

\begin{prop}
{\hspace{-.2cm}.\hspace{.2cm}} \label{prop6}
Let $\Lambda=\left(\begin{array}{ccccc}a_1&,\, \dots\, ,&a_i\ ,\ a_{i+1}&,\, \dots\, ,&a_n\\[.3em]
\tilde{a}_1&,\, \dots\, ,&\tilde{a}_i\ ,\ \tilde{a}_{i+1}&,\, \dots\, ,&\tilde{a}_n\end{array}\right)\in 
{\mathrm{Spec}}\left(\begin{array}{ccc}j_1&,\, \dots\, ,&j_n\\[.5em]
\tj_1&,\, \dots\, ,&\tj_n\end{array}\right)$ and let $e_{\Lambda}$ be a corresponding vector. Then
\begin{itemize}
\item[(a)] We have $a_i^m=1$ for all $i=1,\dots,n$; if $a_i=a_{i+1}$ then $\tilde{a}_i\neq \tilde{a}_{i+1}$.
\item[(b)] If $a_{i+1}=a_i$ and $\tilde{a}_{i+1}=\tilde{a}_i+\epsilon$, where $\epsilon=\pm 1$, then $s_i(e_{\Lambda})=\epsilon e_{\Lambda}$.
\item[(c)] If $a_{i+1}\neq a_i$ or  $a_{i+1}=a_i\And\tilde{a}_{i+1}\neq \tilde{a}_i\pm1$ then 
\[\Lambda'=\left(\begin{array}{ccccc}a_1&,\, \dots\, ,&a_{i+1}\ ,\ a_i&,\, \dots\, ,&a_n\\[.3em]
\tilde{a}_1&,\, \dots\, ,&\tilde{a}_{i+1}\ ,\ \tilde{a}_i&,\, \dots\, ,&\tilde{a}_n\end{array}\right)\in {\mathrm{Spec}}\left(\begin{array}{ccc}j_1&,\, \dots\, ,&j_n\\[.5em]
\tj_1&,\, \dots\, ,&\tj_n\end{array}\right)\ .\]
Moreover, if $a_{i+1}\neq a_i$ then the vector $s_i(e_{\Lambda})$ corresponds to the string $\Lambda'$, see the matrices (\ref{mat-d2'}) with 
$a=a_i$, $b=a_{i+1}$, $\tilde{a}=\tilde{a}_i$ and $\tilde{b}=\tilde{a}_{i+1}$; if $a_{i+1}=a_i$ and $\tilde{a}_{i+1}\neq \tilde{a}_i\pm1$ then the 
vector $s_i(e_{\Lambda})-\frac{1}{\tilde{a}_{i+1}-\tilde{a}_i}e_{\Lambda}$ corresponds to the string $\Lambda'$, see the matrices (\ref{mat-d2''}) 
with $a=a_i=a_{i+1}$, $\tilde{a}=\tilde{a}_i$ and $\tilde{b}=\tilde{a}_{i+1}$.
\end{itemize}
\end{prop}

\subsection{{\hspace{-0.50cm}.\hspace{0.50cm}}Classical content strings}\label{subsec-cla6}

\begin{defi}
{\hspace{-.2cm}.\hspace{.2cm}} \label{def-cont2}
A classical content string $\left(\begin{array}{ccc}a_1&,\, \dots\, ,&a_n\\[.3em]
\tilde{a}_1&,\, \dots\, ,&\tilde{a}_n\end{array}\right)$ is a string of columns of numbers satisfying the following conditions: 
\begin{itemize}
\item[(1)] $\tilde{a}_1=0$ and $a_i^m=1$ for all $i=1,\dots,n$. 
\item[(2)] For all $j>1$: 
if  $\tilde{a}_j\neq 0$ then there exists $i$, $i<j$, such that $a_i=a_j$ and
$\tilde{a}_i\in\{\tilde{a}_j-1,\tilde{a}_j+1\}$.
\item[(3)] If $a_j=a_k$ and $\tilde{a}_j=\tilde{a}_k$ for
$j,k$, $j<k$, then there exist $i_1,i_2\in [j+1,k-1]$ such that $a_{i_1}=a_{i_2}=a_j=a_k$, 
$\tilde{a}_{i_1}=\tilde{a}_j-1$ and $\tilde{a}_{i_2}=\tilde{a}_j+1$.
\end{itemize}
The set of classical content strings we denote by ${\mathrm{cCont}}_m(n)$.
\end{defi}

Here is the classical analogue of the Proposition 6 in \cite{OPdA}.

\begin{prop}
{\hspace{-.2cm}.\hspace{.2cm}} \label{prop7}
Assume that a string of columns of numbers $\ \left(\begin{array}{ccc}a_1&,\, \dots\, ,&a_n\\[.2em]
\tilde{a}_1&,\, \dots\, ,&\tilde{a}_n\end{array}\right)\ $ belongs to the set  
$\, {\mathrm{Spec}}
\left(\begin{array}{ccc}j_1&,\, \dots\, ,&j_n\\[.5em]
\tj_1&,\, \dots\, ,&\tj_n\end{array}\right)$. Then it belongs to the set ${\mathrm{cCont}}_m(n)$.
\end{prop}

\subsection{{\hspace{-0.50cm}.\hspace{0.50cm}}Classical content of a node in a Young $m$-diagram}\label{subsec-cl-cont}

Recall that a Young $m$-diagram, or $m$-partition, is an $m$-tuple of Young diagrams $\lambda^{(m)}=(\lambda_1,\dots,\lambda_m)$. The length of a Young diagram $\lambda$ 
is the number of nodes of the diagram and is denoted by $|\lambda|$. By definition the length of an $m$-tuple $\lambda^{(m)}=(\lambda_1,\dots,\lambda_m)$ is $|\lambda^{(m)}|:=|\lambda_1|+\dots+|\lambda_m|$. 

Let the length of the $m$-tuple be $n$. We place the numbers $1,\dots,n$ in the nodes of these diagrams in such a way that in every diagram the 
numbers in the nodes are in ascending order along rows and columns in right and down directions. This is a standard Young $m$-tableau of 
shape $\lambda^{(m)}$. 

\vskip .2cm
The classical content of a node in a Young diagram is $(s-r)$ when the node lies in the line $r$ and column $s$. To extend this definition to Young 
$m$-diagrams, we have to specify in which diagram of an $m$-diagram the node lies; thus, the content of a node in a Young $m$-diagram is a 
couple of numbers, the first number specifies the diagram (in which the node lies) in the $m$-diagram and the second number gives the content of the 
node in the specified diagram. To relate this information with the spectra of the Jucys--Murphy elements, fix (arbitrarily) a bijection between
the set $\{1,\dots,m\}$ and the set of distinct $m^{\textrm{th}}$ roots of unity; let $\xi_k$ be the root of unity associated with $k\in \{1,\dots,m\}$ by this bijection. 
We define the classical content of a node which lies  in the line $r$ and column $s$ of the $k^{\textrm{th}}$ diagram of the $m$-diagram to be the column
$\left(\begin{array}{c}\xi_k\\[.2em] s-r\end{array}\right)$.

To a standard Young $m$-tableau with length $n$ we associate a string of columns of numbers
$\left(\begin{array}{ccc}a_1&,\, \dots\, ,&a_n\\[.2em]\tilde{a}_1&,\, \dots\, ,&\tilde{a}_n\end{array}\right)$ where 
$\left(\begin{array}{c}a_i\\[.2em] \tilde{a}_i\end{array}\right)$
is the classical 
content of the node in which the number $i$ is placed in the $m$-tableau. This association provides the bijection stated in the following classical analogue of the Proposition 7 in \cite{OPdA}.

\begin{prop}
{\hspace{-.2cm}.\hspace{.2cm}} \label{prop8}
There is a bijection between the set of standard Young $m$-tableaux of length $n$ and the set ${\mathrm{cCont}}_m(n)$.
\end{prop}

Here is an example of a standard Young $2$-tableau with $m=2$ and $n=10$: $$\left(\begin{array}{l} \fbox{\scriptsize{1}}\fbox{\scriptsize{2}}\fbox{\scriptsize{4}}\\[-0.2em]
\fbox{\scriptsize{6}}\fbox{\scriptsize{9}}\\[-0.2em]
\fbox{\scriptsize{7}}\end{array}\ ,
\ \begin{array}{l} \fbox{\scriptsize{3}}\fbox{\scriptsize{8}}\fbox{\scriptsize{10}}\\[-0.2em]
\fbox{\scriptsize{5}}\end{array}\right)\ .$$
The string, associated to this standard Young $2$-tableau, is:
\[\left(\begin{array}{ccccccccccccccccccc}\xi_1&,&\xi_1&,&\xi_2&,&\xi_1&,&\xi_2&,&\xi_1&,&\xi_1&,&
\xi_2&,&\xi_1&,&\xi_2\\[.3em]
0&,&1&,&0&,&2&,&-1&,&-1&,&-2&,&1&,&0&,&2
         \end{array}\right),\]
where $\{ \xi_1,\xi_2\}$ is the set of distinct square roots of unity.
 
\section{{\hspace{-0.55cm}.\hspace{0.55cm}}Construction of representations}\label{suserep-c}

Here we establish an analogue, in the classical setting, of the construction of representations of $H(m,1,n)$ presented in \cite{OPdA}: we define an algebra structure on a 
tensor product of the algebra $\mathbb{C}G(m,1,n)$ with a free associative algebra
generated by the standard $m$-tableaux corresponding to $m$-partitions of $n$. Then, by 
evaluation (with the help of the simplest one-dimensional representation of $G(m,1,n)$) from the right, we build representations.

\subsection{{\hspace{-0.50cm}.\hspace{0.50cm}}Baxterized elements} 

In Subsection \ref{sec5.2} we define the algebra $\mathfrak{T}$. The defining relations of $\mathfrak{T}$, involving the generators $s_i$, can 
be conveniently written in terms of the so-called Baxterized elements; we leave the check to the reader. 

\vskip .2cm
Define, for any $s_i$ among the generators $s_1,\dots,s_{n-1}$ of $G(m,1,n)$, the Baxterized element $s_i(\alpha,\beta)$ by 
\begin{equation}\label{bax-sig-c}
s_i(\alpha,\beta):=s_i+\frac{1}{\alpha-\beta}\ .
\end{equation} 
The parameters $\alpha$ and $\beta$ are called spectral parameters. 

\begin{prop}
 {\hspace{-.2cm}.\hspace{.2cm}}
 \label{prop-bax-c}
The following relations hold:
\begin{equation}\label{rel-bax-c}
\begin{array}{cc}s_i(\alpha,\beta)s_i(\beta,\alpha)=1-\frac{1}{(\alpha-\beta)^2}\ ,\\[0.7em]
s_i(\alpha,\beta)s_{i+1}(\alpha,\gamma)s_i(\beta,\gamma)=s_{i+1}(\beta,\gamma)s_i(\alpha,\gamma)s_{i+1}(\alpha,\beta)\ ,&\\[0.85em]
s_i(\alpha,\beta)s_j(\gamma,\delta)=s_j(\gamma,\delta)s_i(\alpha,\beta)&\textrm{if\ \  $|i-j|>1$\ .}
\end{array}
\end{equation}
\end{prop}

The following Lemma shows that
the original relations for the generators $s_i$ follow from the relations for the Baxterized elements with fixed values of the spectral parameters.

\begin{lemm}
 {\hspace{-.2cm}.\hspace{.2cm}}
 \label{prop-bax2-c}
 Let $A$ and $B$ be two elements of an arbitrary associative unital algebra ${\cal{A}}$. Denote  $A(\alpha,\beta):=A+\frac{1}{\alpha-\beta}$ and
 $B(\alpha,\beta):=B+\frac{1}{\alpha-\beta}$ where $\alpha$ and $\beta$ are parameters.

\medskip
(i) If $A(\alpha,\beta)A(\beta,\alpha)=1-(\alpha-\beta)^{-2}\ $ for some (arbitrarily) fixed values of the parameters $\alpha$ and $\beta$, such that $\alpha\neq\beta$, then $A^2=1\ $. 

\medskip
(ii) If $A^2=1$, $B^2=1$ and $A(\alpha,\beta)B(\alpha,\gamma)A(\beta,\gamma)=B(\beta,\gamma)A(\alpha,\gamma)B(\alpha,\beta)\ $ for some (arbitrarily) fixed values of the parameters $\alpha$,$\beta$ and $\gamma$, such that $\alpha\neq\beta\neq\gamma\neq\alpha$, then 
$ABA=BAB\ $.

\medskip
(iii) If $A(\alpha,\beta)B(\gamma,\delta)=B(\gamma,\delta)A(\alpha,\beta)\ $ for some (arbitrarily) fixed values of the parameters $\alpha$,$\beta$,$\gamma$ and $\delta$, such that $\alpha\neq\beta$ and $\gamma\neq\delta$, then $AB=BA\ $.
\end{lemm}

\subsection{{\hspace{-0.50cm}.\hspace{0.50cm}}Product of the algebra $\mathbb{C}G(m,1,n)$ with a free associative algebra generated by the standard $m$-tableaux corresponding to 
$m$-partitions of $n$} \label{sec5.2}

Let $\lambda^{(m)}$ be an $m$-partition of length $n$. Consider a set of free generators labeled by standard $m$-tableaux of shape $\lambda^{(m)}$; for a standard $m$-tableau $X_{\lambda^{(m)}}$ we denote by $\mathcal{X}_{\lambda^{(m)}}$ the corresponding free generator. 

\vskip .2cm
For a standard $m$-tableau $X_{\lambda^{(m)}}$, we shall denote the entries of the content column
of the node where $i$ is placed by $\left(\begin{array}{c}p(X_{\lambda^{(m)}}|i)\\[.2em] \cc(X_{\lambda^{(m)}}|i)\end{array}\right)$. 

\vskip .2cm
By definition, for any $m$-tableau $X_{\lambda^{(m)}}$ and any permutation 
$\pi\in S_n$, the $m$-tableau $\xl{\lam}{\pi}$ is obtained from the $m$-tableau $X_{\lambda^{(m)}}$ by applying the permutation $\pi$ to the 
numbers occupying the 
nodes of $X_{\lambda^{(m)}}$; for example $\xl{\lam}{s_i}$ is the $m$-tableau 
obtained from $X_{\lambda^{(m)}}$ by exchanging the positions of the numbers $i$ and $(i+1)$ in the $m$-tableau $X_{\lambda^{(m)}}$. 

\vskip .2cm
{}For a standard $m$-tableau $X_{\lambda^{(m)}}$, the $m$-tableau $\xl{\lam}{\pi}$ is not necessarily standard. 
As for the generators of the free algebra, we denote the generator corresponding to the $m$-tableau $\xl{\lam}{\pi}$ by 
$\ml{\lam}{\pi}$ if the $m$-tableau $\xl{\lam}{\pi}$ is standard. And if the $m$-tableau $\xl{\lam}{\pi}$ is not standard 
then we put $\ml{\lam}{\pi}=0$. 

\begin{prop}
{\hspace{-.2cm}.\hspace{.2cm}}
 \label{prop-rel-c}
The following relations:
\begin{itemize}
 \item if $p(X_{\lambda^{(m)}}|i)\neq p(X_{\lambda^{(m)}}|i+1)$ then
\begin{equation}\label{rel-c}s_i\cdot\mathcal{X}_{\lambda^{(m)}}=\ml{\lam}{s_i}\cdot s_i\ ;\end{equation}
\item if $p(X_{\lambda^{(m)}}|i)=p(X_{\lambda^{(m)}}|i+1)$ then
\begin{equation}\label{rel-c2}
\hspace{-0.2cm}\Bigl(s_i+\frac{1}{\cc(X_{\lambda^{(m)}}|i)-\cc(X_{\lambda^{(m)}}|i+1)}\Bigr)\cdot\mathcal{X}_{\lambda^{(m)}}=\ml{\lam}{s_i}\cdot\Bigl(s_i+\frac{1}{\cc(X_{\lambda^{(m)}}|i+1)-\cc(X_{\lambda^{(m)}}|i)}\Bigr)\end{equation}
\end{itemize}
and 
\begin{equation}\label{rel-c3}
\Bigl(t-p(X_{\lambda^{(m)}}|1)\Bigr)\cdot\mathcal{X}_{\lambda^{(m)}}=0
\end{equation}
are compatible with the relations for the generators $t,s_1,\dots,s_{n-1}$ of the group $G(m,1,n)$.
\end{prop}
We explain the meaning of the word ``compatible" in the formulation of the Proposition.

\vskip .2cm
Let ${\cal{F}}$ be the free associative algebra
generated by $\tilde{t}$, $\tilde{s}_1,\dots,\tilde{s}_{n-1}$. The group ring $\mathbb{C}G(m,1,n)$ is naturally the quotient of ${\cal{F}}$. 

\vskip .2cm
Let $\mathbb{C}[\mathcal{X}]$ be 
the free associative algebra whose generators  
$\mathcal{X}_{\lambda^{(m)}}$  range over all standard $m$-tableaux of shape $\lambda^{(m)}$ for all $m$-partitions $\lambda^{(m)}$ of $n$.

\vskip .2cm
Consider an algebra structure on the space $\mathbb{C}[\mathcal{X}]\otimes {\cal{F}}$ for which: (i) the map $\iota_1:x\mapsto x\otimes 1$, 
$x\in \mathbb{C}[\mathcal{X}]$, is an isomorphism of $\mathbb{C}[\mathcal{X}]$ with its image with respect to $\iota_1$; (ii) the map 
$\iota_2:\phi\mapsto 1\otimes \phi$, $\phi\in {\cal{F}}$, is an isomorphism of ${\cal{F}}$ with its image with respect to $\iota_2$; (iii) 
the formulas (\ref{rel-c})-(\ref{rel-c3}), extended by associativity, provide the rules to rewrite elements of the form $(1\otimes \phi)(x\otimes 1)$, 
where $x\in \mathbb{C}[\mathcal{X}]$ and
$\phi\in {\cal{F}}$, as elements of $\mathbb{C}[\mathcal{X}]\otimes {\cal{F}}$. 

\vskip .2cm
The ``compatibility" means that we have an induced structure of an associative algebra 
on the space $\mathbb{C}[\mathcal{X}]\otimes \mathbb{C}G(m,1,n)$. 
More precisely, if we multiply any defining relation of $\mathbb{C}G(m,1,n)$ (the relation is viewed as an element of the free algebra 
${\cal{F}}$) from the right by a generator $\mathcal{X}_{\lambda^{(m)}}$ (this is a combination of the form ``a relation of $\mathbb{C}G(m,1,n)$ times 
$\mathcal{X}_{\lambda^{(m)}}$") and use the ``instructions" (\ref{rel-c})-(\ref{rel-c3}) to move all appearing $\mathcal{X}$'s to the left (the free generator 
changes but the expression stays always linear in $\mathcal{X}$) then we obtain a linear combination of terms of the 
form ``$\ml{\lam}{\pi}$, $\pi\in S_n$, times a relation of $\mathbb{C}G(m,1,n)$". 

\vskip .2cm
We denote the resulting algebra by ${\mathfrak{T}}$. 

\vskip .2cm
It turns out that the product of any Jucys--Murphy element by a free generator $\mathcal{X}_{\lambda^{(m)}}$ is proportional
to $\mathcal{X}_{\lambda^{(m)}}$; the proportionality coefficients are given by the content columns of nodes.

\begin{lemm}
{\hspace{-.2cm}.\hspace{.2cm}}
 \label{prop-rel2-c}
The relations (\ref{rel-c}) imply the relations:
\begin{equation}\label{rel2-c}
\bigl(j_i-p(X_{\lambda^{(m)}}|i)\bigr)\cdot\mathcal{X}_{\lambda^{(m)}}=0\quad\textrm{for all $i=1,\dots,n$\ ,}
\end{equation}
\begin{equation}\label{rel2-c2}
\bigl(\tj_i-\cc(X_{\lambda^{(m)}}|i)\bigr)\cdot\mathcal{X}_{\lambda^{(m)}}=0\quad\textrm{for all $i=1,\dots,n$\ .}
\end{equation}
\end{lemm}

\subsection{{\hspace{-0.50cm}.\hspace{0.50cm}}Representations}

The Proposition \ref{prop-rel-c} provides an effective tool for the construction of representations of $G(m,1,n)$.

\vskip .2cm
Let $\vert\rangle$ be a ``vacuum" - a basic vector of a one-dimensional $G(m,1,n)$-module; for example, 
$s_i|\rangle=|\rangle$ and $t|\rangle=\xi_1|\rangle$. 
Moving, in the expressions $\phi\mathcal{X}_{\lambda^{(m)}}\vert\rangle$, $\phi\in \mathbb{C}G(m,1,n)$, the elements $\mathcal{X}$'s to the left and using the module structure, we build, due to the compatibility, a representation of $\mathbb{C}G(m,1,n)$ on the vector space $U_{\lambda^{(m)}}$ with the basis $\mathcal{X}_{\lambda^{(m)}}\vert\rangle$. We shall, by a slight abuse of notation, denote the symbol $\mathcal{X}_{\lambda^{(m)}}\vert\rangle$ again
by $\mathcal{X}_{\lambda^{(m)}}$. This procedure leads to the following formulas for the action of the generators $t,s_1,\dots,s_{n-1}$ on the basis vectors $\mathcal{X}_{\lambda^{(m)}}$ of $U_{\lambda^{(m)}}$:

 \vskip .2cm
 $\bullet$ if $p(X_{\lambda^{(m)}}|i)\neq p(X_{\lambda^{(m)}}|i+1)$ then
\begin{equation}\label{rep-c}s_i\ :\ \mathcal{X}_{\lambda^{(m)}}\mapsto\ml{\lam}{s_i}\ ,\end{equation}

$\bullet$ if $p(X_{\lambda^{(m)}}|i)=p(X_{\lambda^{(m)}}|i+1)$ then
\begin{equation}\label{rep-c2}
s_i\, :\, \mathcal{X}_{\lambda^{(m)}}
\mapsto{\displaystyle\frac{1}{\cc(X_{\lambda^{(m)}}|i+1)-\cc(X_{\lambda^{(m)}}|i)}}
\mathcal{X}_{\lambda^{(m)}}+\left(1+
{\displaystyle \frac{1}{\cc(X_{\lambda^{(m)}}|i+1)-\cc(X_{\lambda^{(m)}}|i)}}\right)\ml{\lam}{s_i}\, ,
\end{equation}
and 
\begin{equation}\label{rep-c3}
t\ :\ \mathcal{X}_{\lambda^{(m)}}\mapsto p(X_{\lambda^{(m)}}|1)\mathcal{X}_{\lambda^{(m)}}\ .
\end{equation}
As before, it is assumed here that $\ml{\lam}{s_i}=0$ if $\xl{\lam}{s_i}$ is not a standard $m$-tableau.

\paragraph{Remarks.} \textbf{(a)} The constructed representations do not depend (up to isomorphism) on the value of the generators $s_1,\dots,s_{n-1}$ and $t$ on the vacuum $\vert\rangle$. 

\vskip .2cm
\textbf{(b)} In Appendix to this paper we study the classical intertwining operators $\tilde{u}_{i+1}:=\overline{s}_i\tilde{x}_i-\tilde{x}_i\overline{s}_i\in \mathfrak{A}_{m,n}$, $i=1,\dots,n-1$. The image 
under the map $\pi$, defined in (\ref{degaftcycl}), of the element $\tilde{u}_{i+1}$ is $\pi(\tilde{u}_{i+1})=s_i\tj_i-\tj s_i\in\mathbb{C}G(m,1,n)$, $i=1,\dots,n-1$. The action of  $\pi(\tilde{u}_{i+1})$  in a representation $V_{\lambda^{(m)}}$ is: 
\begin{equation}\label{act-int-c}
\mathcal{X}_{\lambda^{(m)}}\mapsto\left(\cc^{(i)}-\cc^{(i+1)}-\delta_{p^{(i)},p^{(i+1)}}\right)\mathcal{X}_{\lambda^{(m)}}^{s_i}\ ,\end{equation}
where $\cc^{(i)}=\cc(X_{\lambda^{(m)}}|i)$, $p^{(i)}=p(X_{\lambda^{(m)}}|i)$, $i=1,\dots,n$; $\delta_{p,p'}$ 
is the Kronecker symbol.

\subsection{{\hspace{-0.50cm}.\hspace{0.50cm}}Scalar product}

The representations of $G(m,1,n)$ given by formulas (\ref{rep-c})--(\ref{rep-c3}) are analogues of the semi-normal representations of the symmetric group. 
In this Subsection we provide analogues for $G(m,1,n)$ of the orthogonal representations of the symmetric group. 

\vskip .2cm 
Let $\lambda^{(m)}$ be an $m$-partition and let $X_{\lambda^{(m)}}$ and $X'_{\lambda^{(m)}}$ be two different standard $m$-tableaux of shape $\lambda^{(m)}$. For brevity we set $\cc^{(i)}=\cc(X_{\lambda^{(m)}}|i)$ and $p^{(i)}=p(X_{\lambda^{(m)}}|i)$ for all $i=1,\dots,n$. Define the following Hermitian 
scalar product on the vector space $U_{\lambda^{(m)}}$:
\begin{equation}\label{scal-prod1}
\langle \mathcal{X}_{\lambda^{(m)}},\mathcal{X}'_{\lambda^{(m)}}\rangle =0\ ,\end{equation}
\begin{equation}\label{scal-prod2}\langle \mathcal{X}_{\lambda^{(m)}},\mathcal{X}_{\lambda^{(m)}}\rangle =\prod\limits_{j,k\colon j<k\, ,\, p^{(j)}=p^{(k)}\, ,\, \cc^{(j)}\notin\{\cc^{(k)}, \cc^{(k)}\pm 1\}}\frac{\cc^{(j)}-\cc^{(k)}-1}{\cc^{(j)}-\cc^{(k)}}\ .
\end{equation}
The so defined scalar product is positive definite.

\begin{prop}{\hspace{-.2cm}.\hspace{.2cm}}
The Hermitian scalar product (\ref{scal-prod1})--(\ref{scal-prod2}) on $U_{\lambda^{(m)}}$
is invariant under the action of the group $G(m,1,n)$ defined by formulas (\ref{rep-c})--(\ref{rep-c3}).
\end{prop}

As a consequence, the operators for the elements of $G(m,1,n)$ are unitary in the basis $\{\tilde{\mathcal{X}}_{\lambda^{(m)}}\}$ where
\[\tilde{\mathcal{X}}_{\lambda^{(m)}}:=\left(\prod\limits_{j,k\colon j<k\, ,\, p^{(j)}=p^{(k)}\, ,\, \cc^{(j)}\notin\{\cc^{(k)}, \cc^{(k)}\pm 1\}}\left(\frac{\cc^{(j)}-\cc^{(k)}-1}{\cc^{(j)}-\cc^{(k)}}\right)^{\frac{1}{2}}\right)\,\mathcal{X}_{\lambda^{(m)}}\]
for any standard $m$-tableau $X_{\lambda^{(m)}}$ of shape $\lambda^{(m)}$. 

\section{{\hspace{-0.55cm}.\hspace{0.55cm}}Completeness\vspace{.1cm}}\label{subsec-comp-c}

Using the $C$-representations of $G(m,1,n)$, we complete our study of the spectrum of the Jucys--Murphy elements of $\mathbb{C}G(m,1,n)$ and of the representation theory of $G(m,1,n)$. 

\begin{prop}
{\hspace{-.2cm}.\hspace{.2cm}}
 \label{prop-cont-spec-c}
The set ${\mathrm{Spec}}
\left(\begin{array}{ccc}j_1&,\, \dots\, ,&j_n\\[.5em]
\tj_1&,\, \dots\, ,&\tj_n\end{array}\right)$, the set ${\mathrm{cCont}}_m(n)$ and the set of standard $m$-tableaux are in bijection.
\end{prop}

\begin{coro}
{\hspace{-.2cm}.\hspace{.2cm}}
 \label{lem-fin-c}
The spectrum of the classical Jucys--Murphy elements is simple in the representations $V_{\lambda^{(m)}}$ (labeled by the $m$-partitions).
\end{coro}

It means that for two different standard $m$-tableaux (not necessarily of the same shape) the elements of ${\mathrm{Spec}}
\left(\begin{array}{ccc}j_1&,\, \dots\, ,&j_n\\[.5em]
\tj_1&,\, \dots\, ,&\tj_n\end{array}\right)$ associated to them by the Proposition \ref{prop-cont-spec-c} are different.

\vskip .2cm
It remains to verify that we obtain within this approach all irreducible representations of the group $G(m,1,n)$.
As the sum of the squares of the dimensions of the constructed representations equals the order of $G(m,1,n)$ 
(this follows from general results on products of Bratteli diagrams, see e.g. the Appendix B in \cite{OPdA2}), the following Proposition completes the verification.

\begin{prop}
{\hspace{-.2cm}.\hspace{.2cm}}
 \label{prop-fin-c}
The representations $V_{\lambda^{(m)}}$ (labeled by the $m$-partitions) of the group $G(m,1,n)$ constructed in the preceding subsection are 
irreducible and pairwise non-isomorphic.
\end{prop}

In particular, we obtain, from the completeness results, the following conclusions.
\begin{itemize}
\item The branching rules for the chain, with respect to $n$, of the groups $G(m,1,n)$ are free of multiplicities.
\item  The centralizer of the sub-algebra $\mathbb{C}G(m,1,n-1)$ in the algebra $\mathbb{C}G(m,1,n)$ is commutative for each $n=1,2,3,\dots$ 
\item The centralizer of the subalgebra $\mathbb{C}G(m,1,n-1)$ in the algebra $\mathbb{C}G(m,1,n)$ is 
generated by the center of $\mathbb{C}G(m,1,n-1)$ and the Jucys--Murphy elements $j_n$ and $\tj_n$. 
\item The subalgebra generated by the Jucys--Murphy elements $j_1,\dots,j_n,\tj_1,\dots,\tj_n$ of the algebra $\mathbb{C}G(m,1,n)$ is 
maximal commutative.
\end{itemize}

\vskip .2cm 
\paragraph{Remark.} For every standard $m$-tableau $X_{\lambda^{(m)}}$ define the element $\mathfrak{p}_{_{\scriptstyle{X_{\lambda^{(m)}}}}}$ of the 
ring $\mathbb{C}G(m,1,n)$ by the following recursion. The initial condition is $\mathfrak{p}_{\varnothing}=1$. Let $\alpha$ be the
node occupied by the number $n$ in $X_{\lambda^{(m)}}$; denote by $\mu^{(m)}$ the $m$-partition of $n-1$ obtained from $\lambda^{(m)}$ by removing the node $\alpha$. 
Let $X_{\mu^{(m)}}$ be the standard $m$-tableau with the numbers $1,\dots ,n-1$ at the same nodes as in $X_{\lambda^{(m)}}$. Then the recursion 
is given by
\begin{equation}\mathfrak{p}_{_{\scriptstyle{X_{\lambda^{(m)}}}}}:=\mathfrak{p}_{_{\scriptstyle{X_{\mu^{(m)}}}}}\hspace{-0.4cm}\prod_{\beta\colon\!\!\begin{array}{l}\scriptstyle{\beta\in 
{\cal{E}}_+(\mu^{(m)})}\\[-0.2em]\scriptstyle{\cc(\beta)\neq \cc(\alpha)}\end{array}  }\hspace{-0.3cm} \frac{\tj_n-\cc(\beta)}{\cc(\alpha)-\cc(\beta)}\prod_{\beta\colon\!\!\begin{array}{l}\scriptstyle{\beta\in 
{\cal{E}}_+(\mu^{(m)})}\\[-0.2em]\scriptstyle{p(\beta)\neq p(\alpha)}\end{array}  }\hspace{-0.3cm}\frac{j_n-p(\beta)}{p(\alpha)-p(\beta)}\ ,\end{equation}
where $\left(\begin{array}{c}p(\beta)\\[.2em] \cc(\beta)\end{array}\right)$ is the classical content of the node $\beta$ and ${\cal{E}}_+(\mu^{(m)})$ is the set of addable nodes of $\mu^{(m)}$ (a node $\beta$ is called addable for an $m$-partition $\mu^{(m)}$ if the $m$-tuple of sets of nodes obtained from $\mu^{(m)}$ by adding $\beta$ is still an $m$-partition). Due to the completeness results of this Section, the elements $\mathfrak{p}_{_{\scriptstyle{X_{\lambda^{(m)}}}}}$ form a complete set of pairwise orthogonal
primitive idempotents of the algebra $\mathbb{C}G(m,1,n)$.

\vskip .2cm
We have a well-defined homomorphism ${\mathfrak{T}}\to \mathbb{C}G(m,1,n)$ which is identical on the generators 
$t,s_1,\dots,s_{n-1}$ and sends $\mathcal{X}_{\lambda^{(m)}}$ to $\mathfrak{p}_{_{\scriptstyle{X_{\lambda^{(m)}}}}}$ for all standard $m$-tableaux
$X_{\lambda^{(m)}}$.

\section*{Appendix.
\hspace{.2cm}  Classical intertwining operators}\addcontentsline{toc}{section}{Appendix A
$\ $ Classical intertwining operators}

Here we describe the intertwining operators in the degenerate cyclotomic affine Hecke algebra $\mathfrak{A}_{m,n}$; they can be used to investigate the spectrum of the
elements $\tilde{x}_i$ in different representations. We discuss the origin of these intertwining operators in the (non-degenerate) affine Hecke 
algebra. We also re-derive the spectrum of the Jucys--Murphy elements $\tj_i$ from the perturbation theory point of view. The intertwining operators can be introduced \cite{WaWa} in the more general context of the wreath Hecke algebra.  
 
\paragraph{1.} The following Proposition follows from the Propositions \ref{prop6} and \ref{prop7}. 
\begin{prop}\label{cospe-cla}
{\hspace{-.2cm}.\hspace{.2cm}}
We have
\begin{equation}\label{glspej} {\mathrm{Spec}}(\tj_i)\subset [1-i,i-1] \quad\textrm{for all $\ i=1,\dots,n$.}\end{equation}
\end{prop}

To give an alternative proof, in the spirit of \cite{IO}, we introduce the following elements of 
the algebra $\mathfrak{A}_{m,n}$:
\begin{equation}\tilde{u}_{i+1}:=\overline{s}_i\tilde{x}_i-\tilde{x}_i\overline{s}_i=
\overline{s}_i(\tilde{x}_i-\tilde{x}_{i+1})+P_{i+1}\ ,\ \ i=1,\dots,n-1,\end{equation} 
where we denoted $P_{i+1}:=\frac{1}{m}\sum_{p=1}^m \ x_i^p x_{i+1}^{-p}$.

\vskip .2cm
The elements $\tilde{u}_i$ are the ``classical intertwining" operators, they satisfy
\begin{equation}\label{propeu}
\left\{\begin{array}{l}
\tilde{u}_{i+1}x_i=x_{i+1}\tilde{u}_{i+1}\ ,\quad \tilde{u}_{i+1}x_{i+1}=x_{i}\tilde{u}_{i+1}\ ,\ \tilde{u}_{i+1}x_j=x_j\tilde{u}_{i+1}\ \ {\text{for}}\ \ j\neq i,i+1\ ,\\[1em]
\tilde{u}_{i+1}\tilde{x}_i=\tilde{x}_{i+1}\tilde{u}_{i+1}\ ,\quad \tilde{u}_{i+1}\tilde{x}_{i+1}=\tilde{x}_{i}\tilde{u}_{i+1}\ ,\ \tilde{u}_{i+1}\tilde{x}_j=
\tilde{x}_j\tilde{u}_{i+1}\ \ {\text{for}}\ \ j\neq i,i+1\ .
\end{array}\right.
\end{equation}
Next, the elements $\tilde{u}_i$ satisfy the Artin relations:
\begin{equation}\label{ybftu}\tilde{u}_i\tilde{u}_{i+1}\tilde{u}_i=\tilde{u}_{i+1}\tilde{u}_i\tilde{u}_{i+1}\ .\end{equation}
One more property of the elements $\tilde{u}_i$ is
\begin{equation}\label{tukva}\tilde{u}_{i+1}^2=-(\tilde{x}_{i}-\tilde{x}_{i+1})^2+P_{i+1}=-\Bigl(\tilde{x}_{i}-\tilde{x}_{i+1}+P_{i+1}\Bigr)
\Bigl(\tilde{x}_{i}-\tilde{x}_{i+1}-P_{i+1}\Bigr) \ .\end{equation}

Therefore, for a polynomial $\chi$ in one variable, we have
\begin{equation}\label{recretx}\tilde{u}_{i+1}\chi (\tilde{x}_i)\tilde{u}_{i+1}=\chi (\tilde{x}_{i+1})\tilde{u}_{i+1}^2
=-\chi (\tilde{x}_{i+1})\Bigl(\tilde{x}_{i}-\tilde{x}_{i+1}+P_{i+1}\Bigr)\Bigl(\tilde{x}_{i}-\tilde{x}_{i+1}-P_{i+1}\Bigr)\ .\end{equation}
The elements $\tilde{x}_i$, $\tilde{x}_{i+1}$ and $P_{i+1}$ commute. In a representation $\rho$, the spectrum of the operator $\rho (P_{i+1})$ is 
contained in $\{ 0,1\}$; taking for $\chi$ the characteristic equation for $\rho (\tilde{x}_i)$, we conclude that
\begin{equation}\label{recretx2}{\mathrm{Spec}}(\rho(\tilde{x}_{i+1}))\subset {\mathrm{Spec}}(\rho(\tilde{x}_{i}))\cup
\Bigl({\mathrm{Spec}}(\rho(\tilde{x}_{i}))+1\Bigr)\cup \Bigl({\mathrm{Spec}}(\rho(\tilde{x}_{i}))-1\Bigr)\ .\end{equation}
Realizing $\tilde{x}_i$ by $\tj_i$ in a representation of the group $G(m,1,n)$ and taking into account the ``initial condition"
$\tj_1=0$ we rederive (\ref{glspej}).

\paragraph{Remark.} The usual degenerate affine Hecke algebra (it corresponds to $m=1$) is distinguished in the sense that the idempotents $P_i$
become trivial: $P_i=1$, in contrast to the degenerate cyclotomic affine Hecke algebras with $m>1$.

\paragraph{2.} Let $\hat{H}_n$ be the affine Hecke algebra, that is, the algebra with the generators $\tau,\sigma_1,\dots,\sigma_{n-1}$
and the defining relations given by the first five lines in (\ref{pres-hec-cyc}). 
The algebra $H(m,1,n)$ is the quotient of $\hat{H}_n$ by $(\tau-v_1)\dots(\tau-v_m)$.
The Jucys--Murphy operators $y_i$ of the algebra $\hat{H}_n$ are defined by $y_1:=\tau$ and, for $i>0$, $y_{i+1}:=\sigma_i y_i \sigma_i$. 

\vskip .2cm
General intertwining operators ${\mathfrak{U}}_{i+1}$, 
$i=1,\dots,n-1$, of the affine Hecke algebra $\hat{H}_n$ are defined to be operators which verify, for all $i=1,\dots,n-1$,
\begin{equation}\label{deinop}
{\mathfrak{U}}_{i+1}y_i=y_{i+1}{\mathfrak{U}}_{i+1}\,,\quad{\mathfrak{U}}_{i+1}y_{i+1}=y_{i}{\mathfrak{U}}_{i+1}\,,
\quad{\mathfrak{U}}_{i+1}y_k=y_{k}{\mathfrak{U}}_{i+1}\ \ \text{for}\ \,k\neq i,i+1.
\end{equation}

\vskip .2cm
The intertwining operators (the solutions ${\mathfrak{U}}_{i+1}:=U_{i+1}$ of (\ref{deinop})) used in \cite{IO} are 
\begin{equation}\label{qinop}U_{i+1}:=\sigma_i y_i-y_i\sigma_i\ .\end{equation}
The operators $U_{i+1}$ satisfy, in addition to (\ref{deinop}), to the
Artin relation,
\begin{equation}\label{ybinop}U_iU_{i+1}U_i=U_{i+1}U_iU_{i+1}\ ,\end{equation} 
and square to the following function of the Jucys--Murphy elements:
\begin{equation}\label{sqinop}U_i^2=-(y_{i+1}-q^2y_i)(y_{i+1}-q^{-2}y_i)\ .\end{equation}    
In contrast to (\ref{tukva}), the right hand side of (\ref{sqinop}) does not contain anything analogous to the projector $P_{i+1}$. We shall explain the
appearance of the projectors in the classical limit. 

\vskip .2cm
The formulas (\ref{JM-Bn1})-(\ref{JM-Bn2}) show that the Taylor series decompositions of the Jucys--Murphy operators of the cyclotomic algebra 
$H(m,1,n)$ begin
 as follows
\begin{equation}\label{deojumc}J_i=j_i+j_i\tj_i\alpha+O(\alpha^2)\ ;\end{equation}
here $\alpha$ is the deformation parameter, $q^2=1+\alpha+O(\alpha^2)$. We ``lift" the formula (\ref{deojumc}) to the affine Hecke algebra by assuming 
that  the Taylor series decompositions of the Jucys--Murphy operators of the affine algebra $\hat{H}_n$ begin 
as follows
\begin{equation}\label{deojuma}y_i=x_i+x_i\tilde{x}_i\alpha+O(\alpha^2)\ ,\end{equation}
where $x_i$ and $\tilde{x}_i$ belong to the degenerate cyclotomic affine Hecke algebra $\mathfrak{A}_{m,n}$, see (\ref{deg-affHec-na})-(\ref{deg-affHec-nc}).

\vskip .2cm
To perform the classical limit one takes into account the following beginning of the Taylor series decomposition of the elements $\sigma_i$:
\begin{equation}\label{deosi}\sigma_i=\overline{s}_i
+\frac{\alpha}{2}+O(\alpha^2)\ .\end{equation}

Note that the operators $U_{i+1}$, given by (\ref{qinop}),
tend, by (\ref{deojuma}) and (\ref{deosi}), to the operators
\begin{equation}u_{i+1}:=\overline{s}_ix_i-x_i\overline{s}_i\equiv \overline{s}_i(x_i-x_{i+1})\ .\end{equation}
The operators $u_{i+1}$ satisfy all the relations for the operators $\tilde{u}_{i+1}$ listed in (\ref{propeu}); but the intertwining operators $u_{i+1}$ do 
not help to understand the spectrum of the images of the elements $\tilde{x}_i$ in some representation.

\vskip .2cm
As noted in  \cite{IO}, the operators ${\mathfrak{U}}_{i+1}:=U_{i+1}f(y_i,y_{i+1})$, where $f$ is an arbitrary function, are intertwining operators which satisfy the Artin relation.
One shows by induction that for any positive integer $L$,
\begin{equation}\label{poyinop}\sigma_i y_i^L-y_i^L \sigma_i=U_{i+1}\cdot\sum_{b=0}^{L-1}\ y_{i\phantom{1}}^{\phantom{L-1}\hspace{-.5cm} b}\, 
y_{i+1}^{\, L-1-b}\ .
\end{equation}
Therefore, the operators $\sigma_i y_i^L-y_i^L \sigma_i$ are intertwining operators for any positive integer $L$.

\vskip .2cm
Under the assumption (\ref{deojuma}), the operators $\sigma_i y_i^m-y_i^m \sigma_i\equiv \sigma_i (y_i^m-1)-(y_i^m-1) \sigma_i$ tend to 0 when 
$\alpha$ tends to 0. These operators are of order $O(\alpha)$. Denote
\begin{equation}\label{cyinop}\tilde{{\cal{U}}}_{i+1}:=\frac{1}{m}\Bigl(\sigma_i\ \frac{y_i^m-1}{\alpha}-\frac{y_i^m-1}{\alpha}\ \sigma_i\Bigr)\ .\end{equation}
The elements $\tilde{{\cal{U}}}_{i+1}$ tend to $\tilde{u}_{i+1}$ when $\alpha$ tends to 0, and the following Lemma allows to recover the result (\ref{tukva}) from the 
perturbative point of view.
\begin{lemm}{\hspace{-.2cm}.\hspace{.2cm}}
We have
\[\tilde{{\cal{U}}}_{i+1}^2=-\Bigl(\tilde{x}_{i}-\tilde{x}_{i+1}+P_{i+1}\Bigr)
\Bigl(\tilde{x}_{i}-\tilde{x}_{i+1}-P_{i+1}\Bigr)+O(\alpha)\ .\]
\end{lemm}

\paragraph{3.} The elements $j_i$ verify $j_i^m=1$, the characteristic equations for the elements $j_i$ are not significant on the classical level. 
It is easy to obtain the characteristic equation for $\tj_i$ starting from a characteristic equation for $J_i$. 
Let $A_0$ be a semi-simple operator on a vector space $V$. Consider a perturbation of $A_0$ of the form 
\begin{equation}\label{poa0}A=A_0+A_0A_1\alpha+O(\alpha^2)\ ,\end{equation}
where $A_1$ is also semi-simple and the operators $A_0$ and $A_1$ commute. Let ${\mathfrak{r}}$ be an eigenvalue of $A_0$ and 
$V_{\mathfrak{r}}$ the corresponding eigenspace. 
The operator $A(\alpha)$ on the space $V_{\mathfrak{r}}$ has, up to the order $\alpha^2$, the form ${\mathfrak{r}}\ \text{Id}+{\mathfrak{r}}A_1$,
and its eigenvalues are ${\mathfrak{r}}+{\mathfrak{r}}{\mathfrak{s}}_l\alpha$ where $\{ {\mathfrak{s}}_l\}$ is the set of eigenvalues of the restriction of $A_1$ 
to $V_{\mathfrak{r}}$.

\vskip .2cm
In the particular situation with $A_0=j_i$, $A_1=\tj_i$ and $A=J_i$, the spectrum of $A$ is, in general, a subset of 
$\left\{ v_lq^{2\eta},\ l=1,\dots,m\ \text{and}\ \eta\in [1-i,i-1]\right\}$ 
\cite{OPdA}
. We first take the limit $v_l\to\xi_l$, $l=1,\dots,m$. Then 
$\xi_l q^{2\eta}= \xi_l+\xi_l \eta\alpha+O(\alpha^2)$ (since $q^2=1+\alpha+O(\alpha^2)$) thus the spectrum of the operator 
$\tj_i$ is a subset of $[1-i,i-1]$ and we recover the Proposition \ref{cospe-cla} from a perturbative point of view.

\end{document}